\documentstyle[12pt,emlines]{article}
\newtheorem{theorem}{Theorem}
\newtheorem{lemma}{Lemma}
\newtheorem{proposition}{Proposition}
\begin{document}
\title{Vector bundles and Arithmetical Groups I.The higher Bruhat-Tits tree
\footnote
{Steklov  Mathematical  Institute,  Russian Academy of Sciences, Vavilova
str.,    42,    117966    Moscow    GSP-1,    Russia;~e-mail     address:
an@parshin.mian.su}
}
\author{A. N. Parshin}
\date{}
\maketitle
Recently   author  has  proposed  a  generalization  of  the  Bruhat-Tits
     buildings for the $n$-dimensional  local  fields  \cite{P}.  It  was
     shown also that
there
     is  a  connection of this construction with classification of vector
     bundles on algebraic surfaces. In this paper we give the proofs of a
     part of results from \cite{P}, concerning with the  construction  of
     Bruhat-Tits  building  for  the  group $PGL(2)$ over two-dimensional
     local field. Some results will be given in a  more  generality,  for
     the  groups $PGL(n)$ or for local fields of arbitrary dimension. The
     applications to vector bundles will be considered in another  paper.
     We refer to \cite{P} for the detailed motivation of our construction
     and   we   restrict   ourselves   only  by  short  remarks  in  this
     introduction.

Usually local fields (fields of dimension 1 in this language) appear from
formal neighbourhoods of the points on  an  algebraic  curve  (or  on  an
arithmetical  curve).  This  can  be  generalized  to  higher dimensional
schemes where the points  will  be  replaced  by  chains  of  irreducible
subvarieties  having  strictly decreasing codimension. And there is a way
to put all them together in a single group called the {\em  adele  group}
of the variety (or the scheme). Many parts of the classical adelic theory
can  be  generalized  to this situation. In particularly this is done for
the cohomology theory of coherent sheaves and for class field theory (see
\cite{FP} for a survey of the existing theory).
     
Here we would like to apply these notions to the theory of buildings. The
main object of the  theory  is  a  simplicial  complex  attached  to  any
reductive  algebraic  group  $G$  defined over a field $K$. There are two
parallel theories for the cases when  $K$  has  no  additional  structure
(being  a  local  field of dimension 0) and when $K$ is a local field (of
dimension 1). They are known as the {\em spherical} and  {\em  euclidean}
buildings  correspondingly  (see  \cite{BT1, BT2} for original papers and
\cite{B, R, T} for the surveys and text-books. There exists also a system
of axioms which cover both of the cases and which is known as  theory  of
$BN$-pairs or the Tits systems
\cite{B}. As was shown in \cite{P} these two theories (for the groups
  $PGL$) are the special cases
     of  the  general construction for the groups over local fields of an
     arbitrary dimension.

The simplest case which still exhibits the main features of the theory is
the case when $G$ is of rank 1 and more precisely is a group $PGL(V)$  of
projective  linear  transformations  of a vector space $V$ of dimension 2
over a field $K$. Here we restrict ourselves by this case because it is
     quite sufficient for the study of vector bundles (of  rank  2)  over
     algebraic  surfaces  (  and supplies also a necessary background for
     the corresponding theory on arithmetical surfaces).

The paper contains four sections. In the first
we remind the main notions and results on local fields  of  dimension  2.
     All  the  algebraic constructions, like the Weyl group or the Bruhat
     decomposition are collected in section 2. The known results for  the
     Bruhat-Tits  tree  are  shortly  discussed  in  section  3. The last
     section contains the main results of the paper --  construction  and
     the
properties of the Bruhat-Tits tree over two-dimensional local field.

A  substantial  part  of  this  work  was  done  during  my visits to the
"Sonderfoschungsbereich  170"  of  the  Mathematical  Institute  of   the
G\"ottingen University. I am very much grateful to the members of the SFB
for  the  possibility  to work there and my special thanks to Hans Opolka
and Samuel Patterson for the hospitality.

Author has tried to generalize the Bruhat-Tits theory starting from the middle
     80's. The talks with L.  Breen,  B.  Seifert,  J.  Tits  and  T.  A.
     Springer were very useful for me at the beginning of this work.
     
     I am greatly indebted to T. Fimmel who taught me how to work with \TeX.

\section{Local Fields}

We begin with the main definition of the higher adelic theory.
\par\smallskip
{\sc Definition 1}.
Let  $K$  and  $k$  be  fields. We say that $K$ is a $n$-{\em dimensional
local field} with $k$ as {\em last residue field} if the  field  $K$  has
the following structure. Either $n = 0$ or $K$ is the quotient field of a
(complete)  discrete valuation ring ${\cal O}_{K}$ whose residue field is
a local field of dimension
$n-1$ with last residue field $k$. If $K', K''$, ... are the intermediate
residue fields from the definition then we will write $K/K'/K''.../k$ for
the structure.The {\em first residue field} will  be  denoted  mostly  as
$\bar K$.

\par\smallskip
In  the sequel we restrict ourselves by the case of $n = 2$ ( the general
case was considered in \cite{P}).
     
A typical example (which is quite sufficient  if  we  have  in  mind  the
applications to algebraic surfaces) is the field of iterated power series
$$ K = k((u))((t)) $$
with an obvious inductive local structure on it
$$ {\cal O}_{K} = k((u))[[t]],\,  \bar K = k((u)). $$
(see  \cite[ch.2]{FP})  for other examples and classification theorem for
complete local fields). Let us mention that  the  choice  of  {\em  local
parameters}  $t,  u  $  in  our  example  does  not follow from the local
structure.

For technical reasons we do {\em not} assume as usual that  the  discrete
valuation  rings which enter in our definition are the complete rings. We
have reduction map $p :{\cal O}_{K} \rightarrow {\bar K} $ and we  denote
by $\wp$ and $m$ the maximal ideals of the local rings
${\cal O}_{K}$ and ${\cal O}_{\bar K}$ correspondingly.  Also we denote
     by $t, u$ the generators of these ideals.
Let
$$ {\cal O}'_{K} = p^{-1}({\cal O}_{\bar K}) $$
be a subring in $K$.

Then  we  have a tower of valuation rings for the valuations $\nu^{(i)}$~
of rank $i = 0, 1, 2$:
$$ {\cal O}_{(0)} \supset {\cal O}_{(1)} \supset {\cal O}_{(2)},$$
where $ {\cal O}_{(0)} = K,~{\cal O}_{(1)} = {\cal O}_{K},~{\cal O}_{(2)}
= {\cal O}'$.

For valuation groups
$$ \Gamma_{K}^{(i)} = K^{*}/({\cal O}_{(i)})^{*},~i = 1,2,
~\Gamma_{K} = \Gamma_{K}^{(2)}  $$
there is a filtration
$$ \Gamma_{K}^{(2)} \rightarrow \Gamma_{K}^{(1)} , $$
which will be reduced to one homomorphism in our case. We  denote  it  by
$\pi$.

This  filtration  defines  on  $\Gamma_{K}$  a structure of {\em ordered}
group.  If  we  need  to  show  the  local  structure   we   will   write
$\Gamma_{K/.../k}$  instead  of $\Gamma_{K}$.If we chose local parameters
$t, u$ of the field $K$ then the order becomes the lexicographical order.
Inside the group $\Gamma_{K}$ there is a subset
$\Gamma_{K}^{+}$ of non-negative elements.

In our situation we have two valuations (of ranks 1 and 2). They will  be
     denoted by $ \nu $ and $ \nu'$ correspondingly.
     
If $K \supset {\cal O}$ is a fraction field of a subring ${\cal O}$ we call
${\cal O}$-submodules $a \subset K$ fractional ${\cal O}$-ideals (or
simply fractional ideals).

\begin{theorem}.
The  local  rings  ${\cal  O}_{(i)}  $  $i  = 0, 1, 2$ have the following
     properties:
\par\smallskip
     
i) $$ {\cal O}'/m  =  k,~K^{*}  =  \{ t \} \{ u \}({\cal O}')^{*},
   ~({\cal O}')^{*}  =  k^{*}(1 +  m);$$

ii) every finitely  generated  fractional  ${\cal  O'}$-ideal  $a$  is  a
principal one and
$$a = m_{i, n} = (u^{i}t^{n}),~i,n \in \mbox{\bf Z}; $$

iii) every infinitely generated fractional ${\cal O'}$-ideal $a$ is equal to
$$a  = \wp_{n} = (u^{i}t^{n} \mid \mbox{for all~}
i \in \mbox{\bf Z}),~n \in \mbox{\bf Z}; $$

iv) if $\wp_{(i,j)} = \wp \mbox{ for (i,j) = (2,1) and } \wp_{(i,j)} = (0)
\mbox{ for (i,j) = (i,0)},$
      $$ {\cal  O}_{(i)} \supset \wp_{(i,i-1)} \supset ... \supset \wp_{(i,1)}
        \supset \wp_{(i,0)}, $$
then
$$ \mbox{Hom}_{{\cal O}'}({\cal O}_{(i)}, {\cal O}_{(j)}) = \left\{
\begin{array}{ll}
 {\cal O}_{(j)}, & i \geq j,    \\
 \wp_{(j,i)}, & i < j.
\end{array}  \right.
$$
\end{theorem}

 {\sc Proof}.
The  multiplicative structure of the field $K$ can be deduced immediately
from the corresponding results for the fields of dimension  1  (see,  for
example,
\cite{S1}).  We get also that  $\nu':K^{*}   \rightarrow
\Gamma_{K}$  is a valuation, if we introduce on $\Gamma_{K}$ the lexicographical
order. Thus for any $x, y \in K^{*}$  we have
\begin{equation}
\label{valuation}
     x  =  ay  ~\mbox{with}  ~a \in {\cal O}' \Longleftrightarrow ~\nu'(x)
\geq \nu'(y).
\end{equation}
Let   now   $a  =  (x_{1},...,x_{n})$  be  a  finitely  generated  ${\cal
O}'$-module. If $x \in a$ then $x = \sum a_{i}x_{i},~a_{i}
\in {\cal O}'$.
      From here we see that $min_{x \in a-(0)} \nu'(x)$  ~exists  and  it
can be achieved for some $x_{0} \in a$. This shows that $a = (x_{0})$.
     
If  $b \subset K$ is an infinitely generated module then for some $i$ the
group $b \otimes {\cal O}_{(i)}$ will be a ${\cal  O}_{(i)}$-module  with
one generator. Let
$i$ be the largest index with this property. Then  $\nu^{(i)}$  has it's
minimum  on $b$ and $\min~\nu^{(j)}$ equals to infinity for $j > i$. This
gives our claim.
     
The last property must be checked only for $i <  j$  (  otherwise  it  is
obvious). Explicitly it means that
     $$ \mbox{Hom}_{{\cal O}'}(K, {\cal O}) = \mbox{Hom}_{{\cal O}'}
     (K, {\cal O}') = (0)$$
and
     $$ \mbox{Hom}_{{\cal O}'}({\cal O}, {\cal O}') = \wp .$$
Both the equalities followed from the results already proved on structure
     of ideals in the ring ${\cal O}'$.
     
{\bf Remark 1}.
These  non-noetherian  rings  play an important role in the whole theory.
Usually the higher local fields appear as fields attached to  some  chain
of  the  subschemes of decreasing codimension \cite[ch.4, 7]{FP} and many
structures related with them can be interpreted in  terms  of  simplicial
stucture  on  the partially ordered set of such chains. It seems that the
rings ${\cal O}_{(i)}$ cannot be described in these terms.  It  would  be
interesting   to   extend   the  simplicial  language  (as  described  in
\cite[ch.7]{FP}) to cover these rings also.

\section{BN-pairs}

Let $G = \mbox{SL}(n, K)$ where $K$ is a two-dimensional local field.  We
put
$$
B = \left( \begin{array}{llll}
            {\cal O}' & {\cal O}'& \dots & {\cal O}'\\
            m         & {\cal O}'& \dots & {\cal O}'\\
                      &          & \dots &           \\
            m         &  m       &  \dots & {\cal O}'
            \end{array} \right),
$$
     
We  denote  in  such  way  the subgroup of $G$ consisting of the matrices
whose entries satisfy  the  written  conditions.  Also  let  $N$  be  the
subgroup of monomial matrices.

\par\smallskip
{\sc Definition 2. }
                   Let

$$ T = B \bigcap N = \left( \begin{array}{lll}
                           ({\cal O}')^{*} & \dots & 0 \\
                                               & \ddots &  \\
                                            0   & \dots & ({\cal O}')^{*}
                                            \end{array} \right)
$$

The group
$$W_{K/{\bar K}/k)} = N/T.$$
will be called the
{\em Weyl group}

We also introduce

$$
P = \left( \begin{array}{lll}
                {\cal O}'& \dots & {\cal O}'\\
                         & \dots &     \\
                         & \dots &  \\
                {\cal O}'& \dots & {\cal O}'
                \end{array} \right) \bigcap G,
$$
$$
A = \left\{ \left( \begin{array} {lll}
            t^{i_{1}} u^{j_{1}} & \dots &   0  \\
            0                                         &  \dots &  0 \\
                                  & \ddots       & \\
                              0  & \dots & t^{i_{n}} u^{j_{n}} \\
                              \end{array} \right),
                \mbox{~for all } k~i_{k},j_{k} \in \mbox{\bf Z} \right\}.
$$
If the matrices in this definition satisfy the additional condition : the
integer   vectors   $(i_{1},   j_{1}),   \dots,   (i_{n},   j_{n})$   are
lexicographically ordered, then we get a subset $A_{+}$. We set also
$$
U = \left( \begin{array}{lll}
              1 & \dots & K \\
                & \ddots &   \\
              0 & \dots & 1
              \end{array} \right)
              $$
\begin{theorem}
We have the following decompositions in the group $G$:
\begin{quotation}
     
          i) the Bruhat decomposition
          $$ G = \bigcup_{w \in W} BwB $$
     
          ii) the Cartan decomposition
          $$ G = \bigcup_{a \in A_{+}} PaP $$
     
          iii) the Iwasawa decomposition
          $$ G = \bigcup_{a \in A} PaU $$
\end{quotation}
          where all the unions are disjoint ones.
     
\end{theorem}
     
{\sc Proof} can be given as a generalization of the known proofs of these
facts for local fields of dimension  1  (  see,  for  example,  \cite[ch.
VI]{G}  for  the  Cartan  and  Iwasawa  decompositions  and \cite[theorem
3.15]{GI} for the Bruhat decomposition). We only outline the  main  steps
here.

{\sc  Existence. } This can be done by standard application of elementary
trasformations to the rows and columns of matrices from  the  group  $G$.
Let $e_{i, j}(\lambda)$ be an elementary matrix with $\lambda$ on
     $(i, j)$-th place. Now let $g = (a_{k, l}) \in G$  and for some
     $~k, i, j$~
$\nu'(a_{k, i}) \le (\mbox{or}<~) \nu'(a_{k, j}) $. Then after a
multiplication from the right by
$e_{i, j}(\lambda)$,~$\lambda = - a_{k, i}^{-1} a_{k, j}$ we get 0 on
 the    $(k, j)$-th place. By~(\ref{valuation}), $\lambda \in {\cal O}'
~(\mbox{or}~m)$. The same fact is true for the multiplication
 by $e_{i, j}(\lambda)$ from the left.
     
     Multiplying the given matrix from $G$ by $e_{i, j}(\lambda)$
     ~with~ $\lambda \in {\cal O}'$~for~$i < j$~and~$\lambda \in m$~for~$i > j$,
     from  the left and from the right we can get a monomial matrix after
    several steps. This gives the Bruhat decomposition.
           In other two cases we need also  to  multiply  by  permutation
 matrices (after a multiplication by a suitable matrix from $T$ they
    belong  to $ SL(n, {\cal O}')$). Also we have to change $m$~on~${\cal
O}'$~ in the second restriction on  matrices  $e_{i,  j}(\lambda)$  given
above.
                   
 {\sc  Uniqueness.}  Let  $L  =  {\cal O}'e_{1} \oplus \dots \oplus {\cal
O}'e_{n}$ be a free ${\cal O}'$-submodule of the space $V$. If
     $x \in V$~ ~$g \in \mbox{GL}(V)$, then we put
   \begin{equation}
     \label{val1}
\nu'(x) = \mbox{min}_{i} \nu'(x_{i}),~\nu'(g) = \mbox{min}_{i, j}
\nu'(a_{i, j}),
     \end{equation}
     where $x = x_{1}e_{1} + \dots + x_{n}e_{n}$.
\par\smallskip
\begin{lemma}. $\nu'(x) \in \Gamma_{K} \cup \infty $ ~and we have:
     
     i) $\nu'(x) = \mbox{min}_{x \in \lambda L}~ \nu'(\lambda),$
     
     ii) $\nu'(g) = \nu'(pgq), ~\mbox{if}~p, q \in
\mbox{Stab}(L) \cong \mbox{GL}(n, {\cal O}'),$
     
     iii) $\nu'(g) = \mbox{min}_{\nu'(x) = 0}~ \nu'(g(x))
        = \mbox{min}_{x \in L}~ \nu' (g(x)).$
\end{lemma}
     \par\smallskip
     These properties can be checked precisely as in the case of discrete
valuation rings. We denote $\nu'(x)$ by $\nu'_{L}(x)$ because it
     depends only on the submodule $L$.

We show how to get the uniqueness for the Bruhat decomposition.

     Let $L_{k} = me_{1} \oplus \dots \oplus me_{k} \oplus {\cal  O}'e_{k
+ 1} \dots \oplus {\cal O}'e_{n},~k = 0, ..., n - 1$ --- free
${\cal O}'$-submodules in $V$. We put
     \begin{equation}
     \label{val2}
      \delta_{rkl}(g) = \mbox{min}_{x \in \wedge^{r} L_{k}}
     \nu'_{\wedge^{r}L_{l}}(\wedge^{r}g(x)),
     \end{equation}
     where  $r  =  1,  ..., n$~and~$\wedge^{r}L_{k}$ is a $r$-th external
power of the module $L_{k}$~in~$\wedge^{r}V$. Then we can prove:
     $$\delta_{rkl}(g') = \delta_{rkl}(g),~\mbox{if}~g' \in BgB,$$
     (since ~$\forall k, B(L_{k}) = L_{k}$) and
     $$\mbox{if}~w, w' \in N~\mbox{and}~\forall r, k, l
~\delta_{rkl}(w) = \delta_{rkl}(w'),~\mbox{then}~w'w^{-1} \in T.$$
     
     It gives our claim. The uniqueness for the Cartan decomposition  can
be proved along the same lines (with one module $L$ instead of all
$L_{k}$) . The uniqueness for the Iwasawa decomposition is a direct
     computation.
     
The proof is finished.
\par\smallskip
     {\bf Remark 2.} The same type decompositions also exist for the group
$\mbox{GL}(n, K)$.
     
     We  also  conjecture that the decompositions of the theorem (and the
known decompositions  for  the  parabolic(parahoric)  subgroups  in  Tits
theory
\cite[ch. IV, \S 2.5]{B}) can be  included in some general theorem
formulated in an appropriate simplicial language using the rings
 ${\cal O}_{(i)}$.
     \par\smallskip
        Let us study the Weyl group $W$ more carefully.
It contains the following elements of order two
$$ s_{i} = \left( \begin{array} {llllllll}
                   1  & \dots  &  0  &    &   & 0 & \dots & 0 \\
                      & \ddots &     &    &   &   &     &    \\
                   0  & \dots  &  1  &    &   &   &     &  0 \\
                   0  & \dots  &     & 0  & 1 &   & \dots  & 0     \\
                   0  &        &     & -1 & 0 &   &  \dots  & 0      \\
                   0  & \dots  &     &    &   & 1 & \dots   & 0  \\
                      &        &     &    &   &   & \ddots &       \\
                   0  &  \dots &  0  &    &   & 0  & \dots  & 1
                   \end{array} \right), i = 1, ..., n - 1;
$$
$$ w_{1} = \left( \begin{array} {rllll}
                   0 & 0 & \dots & 0 & t \\
                   0 & 1 & \dots & 0 & 0 \\
                     &   & \dots &   &    \\
                     &   & \dots &   &   \\
                     &   & \dots & 1 & 0   \\
         -t^{-1} & 0 & \dots & 0 & 0
         \end{array} \right),~w_{2} = \left( \begin{array} {rllll}
                   0 & 0 & \dots & 0 & u \\
                   0 & 1 & \dots & 0 & 0 \\
                     &   & \dots &   &    \\
                     &   & \dots &   &   \\
                     &   & \dots & 1 & 0   \\
         -u^{-1} & 0 & \dots & 0 & 0
         \end{array} \right).
$$

If  $n  =  2$ then we denote by $ w_{0}$ the element $s_{1}$ of the group
$G$. For general $n$ let $S$ be the constructed set of  elements  of  the
Weyl  group.Then  $  \#S  =  n  +  1  $  (  and  $\mbox{rk}(G)  +  m$ for
$m$-dimensional field) and we have

\par\smallskip
\begin{theorem}.
The Weyl group $W$ has the following properties:
\begin{quotation}
     
     i) $W$ is generated by the set $S$ of it's elements of order two,

     ii) there exists an exact sequence
        $$
        0 \rightarrow E(= \mbox{Ker } \Sigma) \rightarrow W_{K/\bar K/k}
       \rightarrow  W_{K}  \rightarrow  1, \\
       $$
        where
       $$ \Sigma: \Gamma_{K} \oplus \stackrel{n} \dots \oplus \Gamma_{K}
       \rightarrow \Gamma_{K} $$
   is a  summation map and $ W_{K}$ is isomorphic to the symmetric group
$\mbox{Symm}_{n}$ of $n$ elements,
     
     iii) the elements $s_{i},~i = 1, \dots, n - 1$ define a splitting of
the exact sequence and the subgroup
 $<s_{1}, \dots, s_{n - 1}>$
acts on $E$ by permutations,
     
     iv) if $n = 2$ then the group $W$ has a presentation
$$W = <w_{0}, w_{1}, w_{2}/w_{0}^{2} = w_{1}^{2} = w_{2}^{2} = e,
          (w_{0}w_{1}w_{2})^{2} = e>,$$
     
     v) the Weyl groups of the group $G$ (for $n = 2$) over the local fields
$K,~K/{\bar  K}, ~K/{\bar  K}/k$ can be related by the following diagram
        $$
        \begin{array}{lllllllll}
  &             & 0               &             & 0            & & & & \\
  &             & \uparrow        &             & \uparrow     &
& & & \\
0 & \rightarrow & \Gamma_{K/\bar K} & \rightarrow & W_{K/\bar K} & \rightarrow
& W_{K} & \rightarrow & 0 \\
 &              & \uparrow        &             & \uparrow     &
& \|    &             &    \\
0 & \rightarrow & \Gamma_{K/\bar K/k}  & \rightarrow & W_{K/\bar K/k}   &
\rightarrow & W_{K} & \rightarrow &  0 \\
 &              & \uparrow        &             & \uparrow     & & & & \\
 &              & \Gamma_{\bar K/k}  & =           & {\bf Z}      & & & & \\
 &              & \uparrow        &             & \uparrow     & & & & \\
 &              & 0               &             & 0              & & & &
 \end{array}
$$

\end{quotation}

\end{theorem}
\par\smallskip
{\sc  Proof}.  The  claims  i)  -  iii)  and  v)  follow  from  a  direct
computation. We have to use the multiplicative structure of the local
     field  $K$  and  the structure of it's valuation rings (theorem 1 of
     the previous section).

     Let us deduce the presentation iv). It  is  easy  to  see  that  the
elements  $w_{0}, w_{1}, w_{2}$ satisfy the conditions of the theorem. In
order to show that there are no other relations we observe that
     according to the claim ii) of the theorem,  the  group  $W$  can  be
presented by some generators $a, b$ (free generators of the subgroup
 $E$), $w_{0}$, and defining relations:
$$ w_{0}^{2} = e,~w_{0}aw_{0} = a^{-1},~w_{0}bw_{0} = b^{-1},~ab = ba. $$
     We  may  assume  that  $a  = w_{0}w_{1}$ and $b = w_{0}w_{2}$. It is
enough to show that these relations are equivalent to  the  relations  of
the claim iv). Indeed, we have

     $$ w_{0}aw_{0} = w_{0}w_{0}w_{1}w_{0} = (w_{0}w_{1})^{-1} $$
and similarly for $b$. Then
     $$ e = (w_{0}w_{1}w_{2})^{2} = w_{0}w_{1}w_{2}w_{0}w_{1}w_{0}w_{0}w_{2} =
     ab^{-1}a^{-1}b, \mbox{i. e.} ~ab = ba. $$
     The  same  formulas  will  also  give  the  equivalence  between our
     defining relations in the opposite direction also.

The theorem is proved.
     
             {\sc Corollary }.
{\em The pair $(W,S)$ is  not a Coxeter group
and furthermore there is no subset $S$ of involutions in  $W$  such  that
$(W,S)$ will be a Coxeter group}.
     
{\sc  Proof.} We prove the second claim at once. Assume that the opposite
     is true and consider the map of $S$ into the quotient-group $W_{K}$.
Choose an involution $s$ from the image of this set. Then the set $S'$ of
elements $S$ mapping to $s$ will generate the Coxeter group $W'$
\cite[ch. IV, \S 1.8]{B}.
By the theorem, it will be an extension of the free abelian group $ E$ of
rank $ > 1$ by a group of order 2. We see simultaneously that there is  a
subgroup  of  $E$  which  has rank $ > 1$ and on which the quotient-group
acts as a multiplication by $ - 1$.
     
     We show that this is impossible for a Coxeter group. Let us consider
the Coxeter diagram of the pair  $(W',S')$(for  its  definition  and  the
properties we need see
\cite[ch. IV, \S 1.9]{B}). It is clear that $\# S' > 1$ and the diagram
contains at least two vertices. If they are not connected by an edge then
     the  group  has  to  contain  a  subgroup  $  \mbox{\bf  Z}/2 \oplus
     \mbox{\bf Z}/2$ ,
which is obviously wrong. If the edge  connecting  the  two  vertices  is
     marked  by some finite number $m > 2$, then our group has to contain
     a
subgroup $\mbox{\bf Z}/m$ which is also impossible. It remains  that  the
diagram  of  our  group  is connected and all the edges are marked by the
symbol $\infty$. Now if we have only two vertices then
     $W'$ cannot include a free abelian subgroup of rank
 $> 1$. If the number of vertices  $ > 2$ then $W'$ must
contain a free subgroup of rank $ > 1$ and this is also impossible.
\par\smallskip
{\bf Remark 3}.
If we consider the Weyl group for the $\mbox{SL}(2, K)$ defined over
$n$-dimensional local field  then it will have $n+1$ generators
$w_{0}, \dots, w_{n}$ and the defining relations will be
$w_{i}^{2} =  1,~(w_{0}w_{i}w_{j})^{2}  =  1$ for all  $i,j$. It is not
a Coxeter group also.
\par\smallskip
{\bf Remark 4}.
Here we see the first basic difference between the Tits theory and  ours.
The formalism of the
$BN$-pairs cannot be applied in our situation, at least, without some
substantial  modifications.  Nevertheless  some  corollaries  of the Tits
axioms are valid, for example the Bruhat  decomposition  (see  theorem  2
above)

In  our  situation there exists still some weaker form of the Tits axioms
(from  \cite[ch.  IV.2]{B}).  More  precisely  they  will  be  true  only
partially and for $n = 2$
 we can replace them by the following formula. Let
$$ w = \left( \begin{array}{ll}
                0 & x^{-1} \\
               - x & 0
                \end{array} \right),~v = v'(x),~w(y) = \left( \begin{array}{ll}
                                              0 & y^{-1} \\
                                             - y & 0
                                              \end{array} \right)
$$
If $s = w_{1}$ then there are three possibilities:
$$
\begin{array}{clll}
v \geq 0       & (BwB)(BsB) & = & BwsB \\
(0,-1)  <  v < 0 & (BwB)(BsB) & = & BwsB \bigcup_{(1,-1) + v \leq v'(y) <
(0,1)  -  v}  Bw(y)B  \\  v  \leq  (0,-1)  &  (BwB)(BsB)  &  =   &   BwsB
\bigcup_{(1,-1) + v \leq v'(y) < (0,1) + v} Bw(y)B
\end{array}
$$
We have the same expression for the diagonal elements $w \in W$. And if
$s = w_{0}, w_{2}$ then the Tits axiom T3
$$ BwBsB \subset BwB \cup BwsB $$

 will be valid. These expressions can be deduced by
     straightforward but rather lengthy computations using the elementary
     transformations from the proof of theorem 2.
     
{\sc  Problem  1.}  To  generalize  the  notion  of $BN$-pair in order to
include  both  the  Tits  axioms  and  the  infinite  decompositions  for
non-Coxeter groups which appear here.

     For the BN-pairs attached to the algebraic groups in the Bruhat-Tits
     theory we also know some finiteness property for the double classes
     $BwB$. Namely, they are the finite unions of the cosets $Bg$.
This property is important for  the definition of the Hecke rings (see
     \cite{IM}). It is easy to see that this property is not preserved in
     the higher dimensions. Thus the usual construction cannot be done in
     our case.
     
     {\sc  Problem  2}.  To  define  an  analog of the Hecke ring for the
     groups over $n$-dimensional local fields for $n > 1$.

\section{ Bruhat-Tits tree over local field of dimension 1}

{\sc The Complex} $\Delta(G, K)$.
First we assume that the field $K$ has no additional structure. Then  the
spherical  building  of $G = PGL(V)$ over $K$ is a complex $\Delta(G, K)$
whose vertices are lines $l$ in the space $V$. All  simplices  of  higher
dimension  are  degenerate and thus it's dimension equals zero. The group
$G(K)$ of rational points over $K$ acts on $\Delta(G, K)$ in a transitive
way.

Let $B$ be a Borel subgroup of $G$,
\begin{eqnarray*}
B
& = &
\left( \begin{array}{ll}
                K^{*} & K \\
                0     & K^{*}
                \end{array} \right) ,
\end{eqnarray*}
Then $B$ is the stabilizer of a line in $V$ and
thus $B$ is a stabilizer of a vertex of $\Delta(G, K)$ . This gives us  a
one to one correspondence between the Borel subgroups and the stabilizers
of the vertices.

The   next   important   object   inside   $\Delta(G,   K)$  is  an  {\em
appartment}~$\Sigma$. To specify it we need to choose a maximal torus $T$
of $G$
\begin{eqnarray*}
T
& = &
\left( \begin{array}{ll}
                K^{*} & 0 \\
                0 & K^{*}
                \end{array} \right)
\end{eqnarray*}

or equivalently a splitting $ V = l_{1} \oplus l_{2} $. This  means  that
the  torus  $T$  will fix the pair of vertices corresponding to the lines
$l_{1}$ and  $l_{2}$.  And  this  pair  is  called  an  appartment.  It's
stabilizer is the normalizer $N$ of the torus $T$,
\begin{displaymath}
N = \left( \begin{array}{ll}
                K^{*} & 0 \\
                0 & K^{*}
                \end{array} \right)
                \bigcup
\left( \begin{array}{ll}
                0 & K^{*} \\
                K^{*} & 0
       \end{array} \right) .
\end{displaymath}

The group $ W = N/T $
is called the {\em Weyl group}. In our case it is of order two and has as
a generator an involution
\begin{eqnarray*}
w_{0}
& = &
\left( \begin{array}{rr}
             0 & 1 \\
            -1 & 0
            \end{array} \right)
\end{eqnarray*}
We see that the appartments are precisely the orbits of the Weyl group W.

     {\sc The Complex} $\Delta(G, K/k)$.
Now we turn to
the case when the field $K$ is a local field (of dimension 1) with residue
field
$k$.
Denote by ${\cal O}$ the
valuation ring ${\cal O}_{K}$, by $u$ a generator of the maximal ideal $m$
and by
$\nu :K \rightarrow \Gamma_{K} \cong \mbox{\bf Z}$ the valuation homomorphism.
Again $G = PGL(V)$.
     
We define the euclidean building
$\Delta(G, K/k)$ as  a one-dimensional complex constructed from classes of
lattices in $V$.
A lattice $L$ is an ${\cal O}_{K}$ -submodule in $V$ which is free and of
rank  2. A class $<L>$ of lattices is the set of all lattices $aL$ for $a
\in K^{*}$. We say that two classes $<L>$ and $<L'>$ are connected as the
vertices by an edge iff for some choice of $L$ and $L'$ we have an  exact
sequence
$$ 0 \rightarrow L' \rightarrow L \rightarrow k \rightarrow 0. $$
This  is  equivalent  to  existence of a maximal totally ordered chain of
${\cal O}$-submodules in $V$ which is invariant  under  multiplication  on
$K^{*}$ and contain  $L$ and $L'$ \cite{BT2}.
Then from the combinatorial point of view $\Delta(G, K/k)$ is a homogeneous
tree
\cite{S2}.
     
We  denote  by $\Delta_{i}(G, K/k)$ the set of $i$-dimensional simplices.
From the construction we deduce easily the following  property  which  we
will use in the sequel:
\begin{quotation}

{\sc Link property}

Let $P \in \Delta_{0}(K/k)$ be represented by a lattice $L$. Then the set
of  edges  going  from $P$ is in one to one canonical correspondence with
the set of lines
${\bf P}(V_{P})$ in the vector space $V_{P} = L/mL$ of dimension 2 over k.
The last
set does not depend on the choice of $L$ (or better to say that there are
canonical isomorphisms between these ${\bf P}(V_{P})$ for different $L$'s
in the same class $<L>$).
\end{quotation}
In  particularly,  if  $k  =  {\mbox{\bf  F}}_{q}$ is a finite field then
$\Delta(G, K/k)$ is locally finite.

The group
$SL(V)$ acts on $\Delta(G, K/k)$ in the following way. It is transitive on the
edges $\in \Delta(G, K/k)$ and has two orbits on the vertices $\in \Delta(K/k)$.

There is a {\em type} of the
vertices $P$ which has two values. To understand this consider  an  exact
sequence
$$ 0 \rightarrow PGL^{+}(V) \rightarrow PGL(V) \rightarrow \mbox{\bf
Z}/2\mbox{\bf Z} \rightarrow 0 $$
where  the  last  homomorphism  is  $\nu(det(.))$  {\em mod} 2. The group
$PGL(V)$ acts on $\Delta_{0}(G, K/k)$ in a transitive way and the group
$PGL^{+}(V)$ has the same two orbits as it's subgroup $SL(V)$.

Now let
\begin{eqnarray*}
B
& = &
\left( \begin{array}{ll}
             {\cal O} & {\cal O} \\
             m & {\cal O}
             \end{array} \right)
\end{eqnarray*}
be the subgroup of $SL(V)$  consisting  of  the  matrices  whose  entries
satisfy  the  written  conditions. Then $B$ is a stabilizer of an edge of
the
$\Delta(G, K/k)$
and all the stabilizers look  like  this  in  an  appropriate  coordinate
system of $V$. The stabilizers of the boundary vertices of the edge are
$$P_{0}  = \left( \begin{array}{rr}
              {\cal O} & {\cal O} \\
              {\cal O} & {\cal O}
       \end{array}       \right),  P_{1}  =  \left( \begin{array}{ll}
             {\cal O} & m^{-1} \\
             m & {\cal O}
     \end{array}   \right)                   $$
We  define  the  subgroup  $N$ as above (so it does not reflect the local
structure on $K$). Instead of the maximal torus we take
\begin{eqnarray*}
T
& = &
\left( \begin{array}{ll}
              {\cal O}^{*} & 0 \\
               0 &  {\cal O}^{*}
               \end{array}  \right)
\end{eqnarray*}
and the Weyl group $W = N/T$ is an extension
$$ 0 \rightarrow \mbox{\bf Z} \rightarrow W \rightarrow \mbox{\bf
Z}/2\mbox{\bf Z} \rightarrow 0 $$
Here we can identify the {\bf Z} with the  valuation  group  $\Gamma_{K}$
and  the  {\bf  Z}/2{\bf  Z} with the previous Weyl group of $G$ over the
field $K$ without local structure.

We have a new involution
\begin{eqnarray*}
w_{1}
& = &
\left( \begin{array}{cl}
              0 & u \\
            -u^{-1}  & 0
       \end{array} \right)
\end{eqnarray*}
and the group $W$ is generated by $w_{0}$ and $w_{1}$.

The appartments $\Sigma$ are now the infinite lines:
\\
\\
\unitlength=1mm
\special{em:linewidth 0.4pt}
\linethickness{0.4pt}
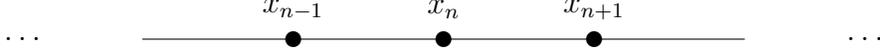
\begin{picture}(126.00,6.00)
\emline{30.00}{2.00}{1}{110.00}{2.00}{2}
\put(90.00,2.00){\circle*{2.00}}
\put(70.00,2.00){\circle*{2.00}}
\put(50.00,2.00){\circle*{2.00}}
\put(50.00,6.00){\makebox(0,0)[cc]{$x_{n-1}$}}
\put(70.00,6.00){\makebox(0,0)[cc]{$x_{n}$}}
\put(90.00,6.00){\makebox(0,0)[cc]{$x_{n+1}$}}
\put(126.00,2.00){\makebox(0,0)[cc]{$\dots$}}
\put(14.00,2.00){\makebox(0,0)[cc]{$\dots$}}
\end{picture}
\\
The group $T$ acts trivially on $\Sigma$ and  it's  stabilizer  coincides
with  $N$.  Thus  appartments  are  the  orbits  of  $W$. The vertices of
$\Sigma$ can be represented by lattices
\begin{equation}
\label{chain}
    x_{n} = <L_{n}>, ~L_{n} = {\cal O} \oplus m^{n}, -\infty < n < \infty
\end{equation}

The action of $W$ on $\Sigma$ can now be  easily  described.  If  $w  \in
\mbox{\bf  Z}$  then  $w$ acts by a translation of even length, and if $w
\not\in \mbox{\bf Z}$ then $w$ acts as an involution with a unique  fixed
point $x_{n_{0}}$:
$$ w(x_{n_{0} + n}) = x_{n_{0} - n}$$

It  can  be  proved that all the appartments look like this and thus they
are in one to one correspondence with the splittings $V  =  l_{1}  \oplus
l_{2}$ of the space $V$.

{\sc  Relations between } $\Delta(G, K)$ {\sc and} $\Delta(G, K/k)$. With
the local field $K$ of dimension 1  we  can  connect  two  local  fields,
namely  $K$ itself and $k$. They are local fields of dimension 0. Thus we
have three buildings attached to $G$ : $\Delta(G, K/k),
\Delta(G, K)$ and $\Delta(G, k)$.

The remark made above (the Link property) shows that the {\em link} of  a
point  $P  \in  \Delta(G,  K/k)$  (  =  the boundary of the $Star(P)$) is
isomorphic to $\Delta(G, k)$. The group
$G(k)$ acts on the last building, $P_{0}$ acts on the link of $P$ and the
isomorphism
between the buildings is an equivariant respective canonical homomorphism from
$P_{0}$ onto $G(k)$ (reduction map {\em mod} $m$).

To formalize the connection with $\Delta(G, K)$ we define a {\em boundary
point} of a tree as a class of half-lines such that intersection  of  any
two half-lines from the class is a half-line in both of them. We have now
an  isomorphism  of  $G(K)$-sets  between  the set of boundary points and
$\Delta(G, K)$. If the half-line is represented by
$L_{n} = {\cal O} \oplus m^{n}, n > 0$ then the corresponding vertex in
$\Delta(K)$ is the line $K \oplus (0)$ in $V$.

It seems reasonable to consider the complexes $\Delta(G, K)$ and
$\Delta(G, K/k)$ together.

Denote by $\Delta_{.}[1](G, K/k)$  the  complex  of  lattices  introduced
above. We define the tree of $G$ as a union
$$ \Delta_{.}(G, K/k) = \Delta_{.}[1](G, K/k) \cup \Delta_{.}[0](G, K) $$
where $\Delta_{.}[1](G, K/k) = \Delta(G, K/k)$ and $\Delta_{.}[0](G, K/k)
     $ is a complex of classes of ${\cal O}$-submodules in $V$ isomorphic
     to
     $K \oplus {\cal O}$
and  we  will  call  the  subcomplex  $\Delta_{.}[0](G,  K/k)$  the  {\em
boundary} of the tree. The definition of the boundary gives a topology on
$\Delta_{0}(G, K/k)$ which is discrete on  both  subsets  $\Delta_{0}[1]$
and $\Delta_{0}[0]$.

Let $P_{n} = <L_{n}>$ and $L_{n} = {\cal O}e_{1} + m^{n}e_{2}$. If $P \in
     \Delta[0]$  is  represented  by  a line $l_{1} = Ke_{1}$ then $P_{n}
\rightarrow P$ since $\cap L_{n} = {\cal O}e_{1}$ belongs to a unique line,
     namely  to  $l_{1}$ (see \cite[ch.II.1.1]{S2}). We can interpret the
     points from $\Delta[0]$ as classes of  ${\cal  O}$-submodules  which
     are  isomorphic  to $K \oplus {\cal O}$ (see lemma 2 below). Then we
     have $P = <L>,~L = Ke_{1} + {\cal O}e_{2}$ instead  of  $l_{1}$  and
     the  definition  of  the  convergence  can  be  given as $ \cup m^{-
     n}L_{n} = L$.

It is easy to
extend it to 1-simplexes. In our case their limits at  infinity  will  be
the  degenerate  simplexes.  Thus  we  have  a  structure of a simplicial
topological space on the tree. It is simply a simplicial  object  in  the
category  of  topological  spaces.  This  topology  is  stronger then the
topology usually introduced to connect these complexes together (see
     \cite{C1}).

Now the connections between the buildings over local fields of  dimension
0 and 1 can be summarized as follows.
$$ \mbox{For any} P \in \Delta_{0}[1](PGL(V), K/k), \mbox{Link}(P) =
\Delta_{.}(PGL(V_{P}), k)       $$
$$   \Delta_{.}[0](PGL(V), K/k) = \Delta_{.}(PGL(V), K) $$
The  last  subset  can be called a star (or link) at infinity. It will be
interesting to define the last notion in  purely  simplicial  terms  (see
remark 6 below).

\section{Bruhat-Tits tree over local field of dimension 2}

 As  above let $G = \mbox{PGL}(V)$ be projective linear group of a vector
space $V$ of dimension 2 over a field $K$ and we assume now that $K$ is a
two-dimensional local field.
     
\par\smallskip
{\sc Definition 3.} {\em The vertices of the Bruhat-Tits tree}.
  $$ \Delta_{0}[2](G, K/{\bar K}/k) = \{ \mbox{classes of}~{\cal O}'
\mbox{-submodules } L \subset
V :  L \cong {\cal O}' \oplus {\cal O}' \}, $$
$$\Delta_{0}[1](G, K/{\bar K}/k) = \{ \mbox{classes of}~{\cal O}'
\mbox{-submodules } L \subset V
:  L \cong {\cal O}' \oplus {\cal O} \}, $$
$$\Delta_{0}[0](G, K/{\bar K}/k) = \{ \mbox{classes of}~{\cal O}'
\mbox{-submodules } L \subset V
: L \cong {\cal O}' \oplus K \}. $$
The two submodules $L$ and $L'$ belong to one class $<L>
= <L'>$, iff $L = aL'$, with $a \in K^{*}$.
     
$$\Delta_{0}(G, K/{\bar K}/k) = \Delta_{0}[2](G, K/{\bar K}/k) \bigcup
\Delta_{0}[1](G, K/{\bar K}/k)
\bigcup \Delta_{0}[0](G, K/{\bar K}/k)$$

We say that the points from $\Delta_{0}[2]$ are the {\em  inner}  points,
the  points  from $\Delta_{0}[1]$ are the {\em inner boundary} points and
the points from $\Delta_{0}[0]$ are the {\em external boundary} points.

Sometimes we will delete $G$ and $K/\bar K/k$ from our notation  if  this
does not lead to a confusion.
     
We  have defined the vertices only. For the simplices of higher dimension
we have the following
\begin{quotation}

{\sc Definition 4.}

Let $\{ L_{\alpha}, \alpha \in I \}$ be a  set of ${\cal O}'$-submodules in
$V$. We say that $\{ L_{\alpha}, \alpha \in I \}$ is a {\em chain} iff:

i) for any $\alpha \in I$ and for any $a \in K^{\star}$ there exists an
$\alpha' \in I$ such that $aL_{\alpha} = L_{\alpha'}$,

ii) the set $\{ L_{\alpha}, \alpha \in I \}$ is totally ordered by the
inclusion.

$\{ L_{\alpha}, \alpha \in I \}$ is a {\em maximal chain} iff it cannot
be included in a strictly larger set satisfying the same conditions i) and ii).

We say that $<L_{0}>, <L_{1}>, ... , <L_{m}>$ belong to a {\em simplex} of
dimension
$m$ iff the $L_{i}, i = 0, 1, ..., m$ belong to a maximal chain of
${\cal O}'$-submodules in $V$. The faces and the degeneracies can be defined in
a standard way (as a deletion or a repetition of some vertex).
\end{quotation}
                                                              
Thus the set $\Delta_{.}(G, K/\bar K/k)$ becomes a {\em simplicial  set}.
The  group  $G = \mbox{PGL}(V)$ acts on ${\cal O}'$-modules. This gives a
simplicial action on $\Delta_{.}(G, K/\bar/k)$.
\par\smallskip
\begin{proposition}.
 The  set  of  all maximal chains of ${\cal O}'$ -submodules in the space
$V$ will be
 exhausted by the following three possibilities:
     \par\smallskip
     i)$~\dots \supset m_{i, n}L \supset m_{i, n}L' \supset m_{i + 1, n}L
     \supset m_{i + 1, n}L' \supset \dots \dots \supset m_{i, n + 1}L
     \supset m_{i, n + 1}L' \supset m_{i + 1, n + 1}L \supset \dots,~i, n
     \in \mbox{\bf Z}, $
     
     where $<L>, <L'> \in \Delta_{0}(G, K/\bar K/k)[2] $ and
     $L \cong {\cal O}' \oplus {\cal O}', ~L' \cong m \oplus {\cal O}'.$
     \par\smallskip
     ii)$~\dots \supset m_{i, n}L \supset m_{i + 1, n}L \supset m_{i + 2, n}
     \supset \dots \dots \supset m_{i, n}L' \supset m_{i + 1, n}L' \supset
  \dots \dots \supset m_{i, n + 1}L \supset m_{i + 1, n + 1}L \supset \dots  ,
     ~i, n \in \mbox{\bf Z},$
     
     where $<L>, <L'> \in \Delta_{0}(G, K/\bar K/k)[1] $ and
     $L \cong {\cal O}' \oplus {\cal O}, ~L' \cong \wp \oplus {\cal O}'.$
     \par\smallskip
     iii) $~\dots \supset m_{i, n}L \supset m_{i + 1, n}L \supset \dots
      \supset m_{i, n + 1}L \supset m_{i + 1, n + 1}L \supset \dots ,
      i, n \in \mbox{\bf Z}$
     
     where $<L>~\in \Delta_{0}(G, K/\bar K/k)[0] $ .
   \end{proposition}
     \par\smallskip
     {\sc Proof.} The chains in the  claim  of  our  proposition  can  be
     completed  by the ${\cal O}'$-modules which are isomorphic to ${\cal
     O}
     \oplus {\cal O}$ and thus do not belong to the modules from
     definition 3. Then a part (with $n = 0$) of the chain of first  type
     will look as follows:
     \begin{equation}
     \label{chain1}
     \dots \supset {\cal O}L \supset \dots \supset L \supset L' \supset
     mL \supset \dots \supset \wp L \supset \dots,
     \end{equation}
     where ${\cal O}L = {\cal O}L' \cong {\cal O} \oplus {\cal O}$~and
     $\wp L = \wp L' \cong \wp \oplus \wp $. There is an isomorphism
     ${\cal O}L / \wp L \cong {\bar K} \oplus {\bar K}.$
     For the same part of the chain of second type we have:
     \begin{eqnarray}
     \dots \supset {\cal O}L \supset \dots \supset L \supset mL \supset
     \dots \supset {\cal O}L' \supset \dots \nonumber \\
     \dots \supset L' \supset mL'\supset
     \dots \supset \wp L \supset \dots,  \label{chain2}
     \end{eqnarray}
     where ${\cal O}L \cong {\cal O} \oplus {\cal O},~{\cal O}L' \cong
     \wp \oplus {\cal O}$~and $\wp L \cong \wp \oplus \wp $. Again there
     exist  isomorphisms
     ${\cal O}L / {\cal O}L' \cong {\bar K},~{\cal O}L' / \wp L
     \cong {\bar K}.$
     The last chain has the following structure:
     \begin{equation}
     \label{chain3}
\dots \supset {\cal O}L \supset \dots \supset L \supset
     mL \supset \dots \supset \wp L \supset \dots ,
     \end{equation}
     where $L \cong K \oplus {\cal O}'$ ~and $ {\cal O}L / \wp L
\cong {\bar K}$.
     
     Let  us  go  to the proof of our proposition. It is easy to see that
     the modules which we have inserted into our chains  are  the  unions
     (intesections)
   of  those  $L_{\alpha}$  which  are  just to the right (left) of them.
    Furthermore,  if  $L_{\alpha'}$  is   the   module,   located   after
    $L_{\alpha}$,
     then  $L_{\alpha}  /  L_{\alpha'}  \cong k $~(theorem 1). It follows
that  all  the  chains  from  i)  --  iii)  are  maximal  ones.  Now  let
$L_{\alpha}$ be an arbitrary maximal chain satisfying to the
     definition 3. We consider three cases:
     
     1)  For  some $\alpha~L_{\alpha} \cong K \oplus {\cal O'}$. Then all
     modules $~m_{i, n}L_{\alpha}$ enter into our chain, i.e. it
will coincide with the chain from iii).
     
     2) For some $\alpha~L_{\alpha} \cong {\cal O}' \oplus {\cal  O}'  $.
     Again all modules $m_{i, n}L_{\alpha}$ belong to the chain, but now
     $L_{\alpha} / mL_{\alpha}$ has dimension 2 over $k$ and since our chain
is  supposed to be a maximal one there exists a module $L'$ between these
two. All $m_{i, n}L'$ are in the chain, which should  coincide  with  the
chain from i).
     
     3)  Now  if  $L_{\alpha}  \cong {\cal O} \oplus {\cal O}'$, then the
   chain contain subchains $\dots \supset m_{i, n}L_{\alpha} \supset
 m_{i + 1, n}L_{\alpha} \supset \dots  $, lying in between the modules
$\wp_{n}L_{\alpha}$ and $\wp_{n + 1}L_{\alpha}$. Choose some $n$.
The intersection of all $m_{i, n}L_{\alpha}$  for varying $i$ gives us a module
 $L'' \supset \wp_{n + 1}L_{\alpha}$. An "empty" place which we have to the
right of $L''$ can be filled out if we
set $L' = \varphi^{-1}({\cal O}_{\bar K})$,
     where $\varphi: L'' \rightarrow L''/\wp L \cong {\bar K}$. Then  all
     "multiples" $m_{i, n}L'$ should be presented in the chain because
of  it's  maximality.  We  see  that the chain constructed in such way is
equal to the chain from ii).
     
     The proposition is proved.
     \par\smallskip
 {\sc Corollary.} {\em The simplicial set $\Delta_{.}$ is a  disconnected
union of it's subsets $\Delta_{.}[m], ~m = 0, 1, 2$. The dimension of the
subset $\Delta_{.}[m]$ equals to 0 for $m = 0$ and 1 for $m = 1, 2$} .
  \par\smallskip
     This  is obvious. We need only note that all vertices of any simplex
   can be represented by the modules of the same type and that in the
case of subset $\Delta_{.}[0]$ the chains of the type iii)  contain  only
{\em one} class of modules.
     \par\smallskip
{\sc Definition 5.}{\em The projection map}.

For any  ${\cal O}'$-module $L$ we have a ${\cal O}$-module
$M = L \otimes _{\cal O'} {\cal O}$ . This gives a map
$$ \pi : \Delta_{.}(K/{\bar K}/k) \rightarrow \Delta_{.}(K/{\bar K}) $$
in  the  tree of the same group $G$ over the field $K$, which we consider
as a local field of dimesnion 1 over ${\bar K}$.
\par\smallskip
\begin{proposition}.
The map $\pi$ has the following properties:
                                          
i) $\pi$ is a simplicial $G$-equivariant surjective map,

ii) $\pi$ induces a bijective map of the  set  $  \Delta_{.}  (G,  K/\bar
K/k)[0]$ onto the set $\Delta_{.}(G, K/\bar K)[0]$,
                                
iii)  if  $\sigma = <L> \in \Delta_{0}(G, K/\bar K)[1]$, then there exist
simplicial and \\
$\mbox{Stab}(<L>)$-equivariant isomorphisms
$$ \pi^{-1}(\sigma) \bigcap \Delta_{.}(G, K/\bar K/k)[2] \cong
\Delta_{.}(\mbox{PGL}(L/\wp L),{\bar K}/k)[1],$$
$$ \pi^{-1}(\sigma) \bigcap \Delta_{.}(G, K/\bar K/k)[1] \cong
\Delta_{.}(\mbox{PGL}(L/\wp L),{\bar K}/k)[0],$$
     where $L/\wp L$ is a vector space of dimension 2 over $\bar K$. Also
we have
$$ \pi^{-1}(\sigma) \bigcap \Delta_{.}(G, K/\bar K/k)[0]
= \emptyset,$$
     
iv) if two vertices from $\Delta_{0}(G, K/\bar K/k)[2]$ are connected  by
an edge then they belong to the same fiber of the map $\pi,$
     
v)  the image of any edge $\sigma \in \Delta_{1}(K/\bar K/k)[1] $ will be
(non-degenerate) edge in $\Delta_{.}(K/\bar K),$

vi) if $\sigma = (\dots \supset L \supset L' \supset \dots ) \in
     \Delta_{1}(G, K/\bar K)[1]$, then   $\pi^{-1}(\sigma)$ consists of
  one edge, connecting vertices from $\pi^{-1}(<L>) \bigcap \Delta_{0}(G,
     K/\bar K/k)[1]$ and  $\pi^{-1}(<L'>)  \bigcap  \Delta_{0}(G,  K/\bar
     K/k)[1]$.
     
\end{proposition}
\par\smallskip
    {\sc Proof.} The property i) is obvious.
  Let  $<L>  \in \Delta_{0}(K/\bar K/k)[0]$~ and let $l \subset L$ be the
set of elements from $L$, divisible in $L$ by all $a \in K^{*}$.
 \par\smallskip
 \begin{lemma}.
      The correspondence $<L> \mapsto l \subset V$~is a bijection between
     $\Delta_{0}[0] $ and $\mbox{\bf P}(V)$.
 \end{lemma}
 \par\smallskip
     {\sc Proof}. If $L = Ke_{1} \oplus {\cal O}'e_{2},$ then
$l = Ke_{1}$   and depends only on class $<L>$ . We can get all the lines
in such way. Now let $L = Ke_{1} \oplus {\cal O}'e_{2},~M = Ke_{1} \oplus
     {\cal O}'e_{2}'$. If $e_{2}' = ae_{1} + be_{2},~a,b \in K$, then  $M
     = Ke_{1} + {\cal O}'be_{2} = bL$,~\\
i.e.~$<L>~=~<M>$.
     
     This  is  also  true  for  $\Delta_{.}(K/\bar  K)[0]$. The claim ii)
follows   since   the   projection   commutes   with   the    constructed
correspondence.
 \par\smallskip
\begin{lemma}.
      Let $L$ be a ${\cal O}'$-submodule in $V$  and
$L \cong {\cal O}' \oplus {\cal O}'$. Then for any point $P \in
\Delta_{0}[2]$~there exists a unique module $L'$ such that $<L'>~= P,\\
     ~L' \subset L$~
     and one of the following equivalent conditions are true:
     
     i) ~$L' \not \subset mL,$
     
     ii)~$L/L' \cong {\cal O}'/a$, where $a$ is a principal ideal,
     
     iii)~$L/L'$ is a module of rank 1.
 \end{lemma}
 \par\smallskip
     {\sc Proof}. Take some module $L''$ in the class of the vertex
     $P$. By the Cartan decomposition (theorem 2)  there exists
    a  basis $e_{1, 2}$ in $V$ such that $L = {\cal O}'e_{1} \oplus {\cal
    O}'e_{2},
     ~L''  =  a_{1}e_{1}  \oplus  a_{2}e_{2},~a_{1, 2}$ are principal
ideals.
     The standard arguments (see \cite[ch. II, \S1.1]{S2}) give the claim
  of the lemma.
     
     We now prove property iii). Fix a module $L_{0} \cong {\cal O}'
     \oplus {\cal O}',~L_{0} \otimes_{{\cal O}'} {\cal O} = L$, i.e.
   ~$<L_{0}> \in
     \pi^{-1}(\sigma)  \cap  \Delta_{0}[2]$. If $P \in \pi^{-1}(\sigma)$,
then according to lemma the vertex $P$ can be  represented  as  $<L'>$  .
Then $L' = m_{i, n}e_{1} +
     {\cal  O}'e_{2},~L_{0}  = {\cal O}'e_{1} + {\cal O}'e_{2}$~ and from
the equality
     $\pi <L_{0}>~ = \pi <L'>~ =~ <L>$~ we get that $n = 0$. It follows
that $L' \supset \wp L_{0}$~and~ $L'$~defines a free
${\cal O}_{\bar K}$-module~$L'/\wp L_{0} \subset L_{0}/\wp L_{0}
     \subset L/\wp L$~of rank 2 in space $L/\wp L$.
     
     This correspondance gives the first bijection from iii). It is  easy
to  see that it preserves the simplicial structure of both sets and it is
equivariant under the stabilizer of the vertex $<L>$.
     
     To construct the second bijection from ii) we take $P \in \pi^{-1}(\sigma)
     \cap \Delta_{0}[1]$. If $P =~ <L'>$ then $<L' \otimes {\cal O}>~ =~ <L>$.
     Changing the module $L'$ to an equivalent one we can assume that
     $L' \otimes {\cal O} = L$ and $L' \subset L$. All such modules $L'$
     can be transformed into one by a multipiliction by some
    $a \in {\cal O}^{*}$. We have a map
     $L' \rightarrow L/ \wp L$. The image $\mbox{Im}~L'$ will be
a ${\cal O}_{\bar K}$
     -module in $L/ \wp L$ isomorphic to ${\bar K} \oplus {\cal  O}_{\bar
     K}$.  As  we  saw  the  class  $<\mbox{Im}~L'>$ will be defined in a
     unique way. It defines  a  point  in  $\Delta_{0}(\mbox{PGL}(L/  \wp
     L),~{\bar  K}/k)[1]$.  The  constructed  correspondence  will  be  a
     bijection with the properties
   we need.
     
     The last claim from iii) follows from the property ii) proved above.

     To get iv) it is enough to apply lemma 2 and proposition 1, i).
     
     The  property  v)  can  be  seen  from  the   description   of   the
  chain~(\ref{chain2}) which represents the edge $\sigma$.
 We need only to take it's quotient by the ideal $\wp$.
     
     We check  now the last property from the proposition. Let $P =~<L>,\\
~Q =~<L'>$
   be  two vertices of the tree $\Delta_{.}(K/\bar K)$ connected by an edge
 $\sigma$. It is posible to choose a basis in $V$ such that $L = {\cal O}e_{1}
     + {\cal O}e_{2},~L' = \wp e_{1} + {\cal O}e_{2}$.
     Then ${\cal O}'$-modules $M = {\cal O}'e_{1} + {\cal O}e_{2}$ and~
     $M' = \wp e_{1} + {\cal O}'e_{2}$ will represent the boundary
points  of the fibers $\pi^{-1}(P)$~and~$\pi^{-1}(Q)$ correspondingly. By
the
  proposition 1, ii) they are connected by an edge which is  mapped  onto
 an edge $\sigma$.
     Thus the set $\pi^{-1}(\sigma)$~is not empty.

It consists of only one edge.
 To  make this clear we denote by $M, M'$ the modules which represent the
vertices of the edge lying over $\sigma$. By the proposition 1, ii)  they
belong  to  a  chain  as in (\ref{chain2}). Now we remark that the set of
lines of the space $L/ \wp L$ is  bijective  to  the  following  sets  of
simplices of our trees:
     
\begin{itemize}
\item  the set of the edges from $\Delta_{.}(K/\bar K)$ going out from the
   vertex $P$~(link property, see section 3).
     
\item   the set $\pi^{-1}(P) \cap \Delta_{1}[1]$~ of the boundary points
  of the fiber $\pi^{-1}(Q)$~(the bijection constructed above).
\end{itemize}
     
From   the  definition  of  the  bijection  we  conclude  that  the  line
corresponding to the vertex $<M>$,
     coincides with the line corresponding to the edge  $\sigma$,  i.  e.
the  vertex  $<M>$ will be defined uniquely. As this is true also for the
vertex $<M'>$ we get that the edge connecting them will be  also  defined
in an unique way.
     
     The proposition is proved.
     \par\smallskip
{\sc  Corollary  1}.{\em  Any  vertex  $P  \in \Delta_{0}[1]$~ belongs to
precisely one edge}.
\par\smallskip
{\sc Corollary 2}. {\em If  $P  \in  \Delta_{0}(K/\bar  K)[1]$  then  the
stabilizer
$G_{P} \subset G$ of the vertex $P$ acts on the fiber $\pi^{-1}(P)$
by the reduction map}
$$ G_{P} \cong \mbox{SL}(2, {\cal O}_{K}) \rightarrow \mbox{SL}(2, {\bar K}).$$

     Here we have fixed the modules $L$ with $<L> = P$ and $L_{0}$ with
$<L_{0}>~\in \pi^{-1}(P) \cap \Delta_{0}[2]$  .
 \par\smallskip
    We  see  that  "inside"  our construction there are five trees coming
from the dimensions $\leq 1$, namely
$$\Delta_{.}(K/{\bar K}),
\Delta_{.}(K), \Delta_{.}({\bar K}/k), \Delta_{.}({\bar K}), \Delta_{.}(k).$$
The first one is the tree which is the target of the projection map
$\pi$, the second one is the external boundary
and the three last trees will occur infinitely many times.

Thus the constructed simplicial set will be a disconnected union of  it's
connected  components.  The  $\Delta_{.}[2]$-piece  of  our  tree  is  an
infinite disconnected union of the usual Bruhat-Tits trees  =  fibers  of
the   map   $\pi$.   The   $~\Delta_{.}[1]$-piece  will  be  an  infinite
disconnected union of the edges.

In order to change this and to have a possibility
to pass from one fiber to another one has to use
some topology which will be a generalization  of  the  topology  we  have
introduced in section 3.

\par\smallskip
{\sc  Definition  6.}We  say  that  a  sequence $P_{n} \in \Delta_{0}[2]$
converges to $Q \in \Delta_{0}[1]$ iff there is a  basis  of  $V$  and  a
sequence  $i(n)$  of  integers  such that $i(n) \rightarrow \infty$ as $n
\rightarrow \infty$ and for large $n$
$P_{n}$ and
$Q$ can be represented by the following modules
$$  P_{n} = <{\cal O}' \oplus m^{i(n)}>, ~Q = <{\cal O} \oplus {\cal O}'>$$
Also a sequence $Q_{n}$ from $\Delta_{0}[1]$ converges to a point $R$ from
$\Delta_{0}[0]$ iff in some basis and for some sequence $i(n)$ as above
$$ Q_{n} = <{\cal O}' \oplus {\wp}^{i(n)}>,~R = <K \oplus {\cal O}'>$$
Combining these two definitions  we  can  get  also  a  condition  for  a
sequence  of  points  from  $\Delta_{0}[2]$  to  converge to a point from
$\Delta_{0}[0]$.

We introduce a topology on $\Delta_{0}(G, K/\bar K/k)$~as a discrete  one
on any of the sets $\Delta_{0}[m]$ and for which the sequences introduced
above
 are the only convergent sequences on the whole set.
$\Delta_{0}$. The convergence on the set of simplices
$\Delta_{1}$
can be defined as convergence of their vertices.
\par\smallskip
  \begin{theorem}. $\Delta_{.}(G, K/\bar K/k)$ is a simplicial
topological  space.  Let  $\mid  \Delta_{.}  \mid  $  be it's geometrical
realization. Then
\begin{quotation}
i) $\mid \Delta_{.} \mid$ is a connected contractible  topological  space
of dimension 1 having a cell structure,

ii)  if  $x  \in  \mid  \Delta_{.}  \mid$  then  $x$  has a neighbourhood
homeomorphic to an interval,  if  $x  \notin  \Delta_{0}[2]$,  and  to  a
bouquet  of  (finite  number  if  \\ $k = {\mbox{\bf F}}_{q}$ ) intervals
otherwise.

iii) the group $G$ acts on $\mid \Delta_{.} \mid$ by homeomorphisms

iv) $\mid \pi \mid$ is a continous map

v) if $\sigma = <L> \in \Delta_{0}(G, K/\bar K)[1]$ then the fiber
$\pi^{-1}(\sigma)$
     is isomorphic to $\Delta_{.}(\mbox{PGL}(L/ \wp L), K/\bar K)$  as  a
simplicial topological space.
\end{quotation}
\end{theorem}
We  refer to \cite{D} for the notions of simplicial topological space and
it's geometrical realization.

{\sc Proof} can be given by a direct check
with an application of the proposition 2 and of the  corresponding  facts
for the trees $\Delta_{.}(K/\bar K)$
~and~$\Delta_{.}(\bar K/k)$ related to the local fields of dimension 1.

{\bf Remark 5}. $\mid \Delta_{.} \mid$ is not a CW-complex even if $n =
1$ and
$k = {\mbox{\bf F}}_{q}$
but it is a closure finite complex. Also we note that
$\mid \Delta_{.}(K/\bar K/k)
\mid $ is not a compact space just as in the  case  of  local  fields  of
dimension 1.

We  can  make  the results proved more transparent by drawing all that in
the following picture where the dots of different  kinds  belong  to  the
different
$\Delta_{.}[m]$-pieces of the tree:
\\
\\
\unitlength=1.00mm
\special{em:linewidth 0.4pt}
\linethickness{0.4pt}
\begin{picture}(127.00,112.00)
\put(15.00,45.00){\circle{4.00}}
\put(125.00,45.00){\circle{4.00}}
\put(34.00,22.00){\circle{4.00}}
\put(65.00,13.00){\circle{4.00}}
\put(94.00,15.00){\circle{4.00}}
\put(113.00,27.00){\circle{4.00}}
\emline{103.00}{45.00}{1}{37.00}{45.00}{2}
\put(48.00,45.00){\circle*{2.00}}
\put(71.00,45.00){\circle*{2.00}}
\put(94.00,45.00){\circle*{2.00}}
\put(62.00,29.00){\circle*{2.00}}
\put(44.00,34.00){\circle*{2.00}}
\emline{44.00}{34.00}{3}{48.00}{45.00}{4}
\emline{71.00}{45.00}{5}{62.00}{29.00}{6}
\put(125.00,110.00){\circle{4.00}}
\put(15.00,110.00){\circle{4.00}}
\put(32.00,110.00){\circle*{2.00}}
\put(38.00,110.00){\circle*{2.00}}
\emline{38.00}{110.00}{7}{32.00}{110.00}{8}
\put(48.00,110.00){\circle{2.00}}
\put(53.00,110.00){\circle{2.00}}
\put(59.00,110.00){\circle{2.00}}
\put(68.00,110.00){\circle*{2.00}}
\put(75.00,110.00){\circle*{2.00}}
\emline{75.00}{110.00}{9}{68.00}{110.00}{10}
\put(85.00,110.00){\circle{2.00}}
\put(93.00,110.00){\circle{2.00}}
\put(103.00,110.00){\circle*{2.00}}
\put(110.00,110.00){\circle*{2.00}}
\emline{110.00}{110.00}{11}{103.00}{110.00}{12}
\emline{97.00}{110.00}{13}{94.00}{110.00}{14}
\emline{92.00}{110.00}{15}{86.00}{110.00}{16}
\emline{84.00}{110.00}{17}{81.00}{110.00}{18}
\emline{62.00}{110.00}{19}{60.00}{110.00}{20}
\emline{58.00}{110.00}{21}{54.00}{110.00}{22}
\emline{52.00}{110.00}{23}{49.00}{110.00}{24}
\emline{47.00}{110.00}{25}{43.00}{110.00}{26}
\put(43.00,103.00){\circle*{2.00}}
\put(52.00,98.00){\circle*{2.00}}
\put(62.00,102.00){\circle*{2.00}}
\put(80.00,104.00){\circle*{2.00}}
\put(89.00,99.00){\circle*{2.00}}
\put(98.00,103.00){\circle*{2.00}}
\put(74.00,97.00){\circle*{2.00}}
\emline{74.00}{97.00}{27}{80.00}{104.00}{28}
\put(68.00,91.00){\circle*{2.00}}
\put(79.00,88.00){\circle*{2.00}}
\put(74.00,91.00){\circle{2.00}}
\put(48.00,88.00){\circle*{0.00}}
\put(48.00,87.00){\circle*{0.00}}
\put(48.00,87.00){\circle*{2.00}}
\emline{48.00}{87.00}{29}{52.00}{98.00}{30}
\put(88.00,84.00){\circle*{2.00}}
\emline{88.00}{84.00}{31}{79.00}{88.00}{32}
\put(70.00,67.00){\circle{4.00}}
\put(116.00,91.00){\circle{4.00}}
\put(98.00,72.00){\circle{4.00}}
\put(40.00,77.00){\circle{4.00}}
\put(119.00,14.00){\makebox(0,0)[cc]{$\Delta_{.}(K/\bar K)$}}
\put(119.00,66.00){\makebox(0,0)[cc]{$\Delta_{.}(K/\bar K/k)$ }}
\put(6.00,78.00){\makebox(0,0)[rc]{$\Big\downarrow \pi$}}
\put(32.00,45.00){\circle*{2.00}}
\emline{37.00}{45.00}{33}{28.00}{45.00}{34}
\put(32.00,45.00){\circle*{2.00}}
\put(81.00,25.00){\circle*{2.00}}
\emline{81.00}{25.00}{35}{62.00}{29.00}{36}
\put(51.00,104.00){\circle{2.00}}
\put(89.00,104.00){\circle{2.00}}
\emline{90.00}{105.00}{37}{93.00}{109.00}{38}
\emline{93.00}{109.00}{39}{93.00}{109.00}{40}
\emline{93.00}{109.00}{41}{93.00}{109.00}{42}
\emline{52.00}{105.00}{43}{53.00}{109.00}{44}
\end{picture}
\\
     \begin{center}
     Pic. 1
     \end{center}

Usually  the  buildings  are  defined as combinatorial complexes having a
system of subcomplexes called appartments (see, for example,
\cite{R, T}).
We show how to introduce them in our case.
\par\smallskip
{\sc Definition 7.} Let us fix a basis
$e_{1},  e_{2} \in V$.  The {\em appartment},~defined by
this basis is the following set
$$ \Sigma_{.} = \bigcup_{0 \leq m \leq 2} \Sigma_{.}[m], $$
where
$$ \Sigma_{0}[m] = \left\{ \begin{array}{l}
<L> \mid L = a_{1}e_{1} \oplus a_{2}e_{2}, \\
\mbox{where }
a_{1}, a_{2}~\mbox{are } {\cal O}'\mbox{-submodules in } K \\
\mbox{and there exists a permutation }
~s ,  \\
\mbox{such that }
a_{s(1)} \cong {\cal O}_{(2)} = {\cal O}',~a_{s(2)} \cong {\cal O}_{(m)}
\end{array} \right\}.
$$
$ \Sigma_{.}[m] $ is the minimal subcomplex having $\Sigma_{0}[m]$ as vertices.

\par\smallskip
Let us denote the edge connecting the vertices $P$~and~ $Q$~by $\sigma(P, Q)$.
\begin{proposition}.
 In some basis we have the following relations:

i)  if
$$
\begin{array}{ccccccc}
x_{i,n} & = & <m_{i,n} \oplus {\cal O}'> & = & <{\cal O}' \oplus m_{-i,-n}>
,& &  \\
     
y_{n} & = & <m_{i,n} \oplus {\cal O}> & = & <m_{j,n} \oplus {\cal O}>
& = & <{\cal O}' \oplus \wp^{-n} >  ,  \\

z_{n}
& = & <{\cal O} \oplus m_{i,-n}> & = & <{\cal O} \oplus m_{j,-n}>
     & = & <\wp^{n} \oplus  {\cal O}' >,
\end{array}
$$
$$ x_{0} =~<K \oplus {\cal O}'>,~x_{\infty} =~<{\cal O}' \oplus K>, $$
then
$$ \Sigma_{0}[2] = \{ x_{i, n} \mid i, n  \in \mbox{\bf Z} \},~
 \Sigma_{1}[2] = \{ \sigma(x_{i, n},~x_{i + 1, n}) \mid i, n
\in \mbox{\bf Z} \},$$
$$ \Sigma_{0}[1] = \{ y_{n}, z_{n} \mid n \in \mbox{\bf Z} \},
~\Sigma_{1}[1] = \{ \sigma(y_{n},~z_{n}) \mid n \in \mbox{\bf Z} \}, $$
$$ \Sigma_{.}[0] = \{ x_{0}, x_{\infty} \}, $$

ii) let $\mbox{Stab}(\sigma)$ be a stabilizer of a  simplex  $\sigma$  in
the subgroup $SL(V)$.
     Then
$$ \mbox{Stab}(x_{i,n}) = \left( \begin{array}{ll}
                           {\cal O}' & m_{i, n}\\
                           m_{-i, -n} & {\cal O}'
                           \end{array} \right),
\mbox{Stab}(\sigma(x_{i, n},~x_{i + 1, n})) = \left( \begin{array}{ll}
                                                 {\cal O}' &  m_{i + 1, n}\\
                                                 m_{-i, -n} &      {\cal O}'
                                                 \end{array} \right) $$
$$\mbox{Stab}(z_{n}) = \left( \begin{array}{ll}
                     {\cal O} & \wp^{n} \\
                       \wp^{-n + 1}    & {\cal O}'
                       \end{array} \right),  \mbox{Stab}(y_{n}) = \left(
\begin{array}{ll}
   {\cal O}'& \wp^{n + 1} \\
           \wp^{-n} &  {\cal O}
     \end{array} \right), $$
$$\mbox{Stab}(\sigma(y_{n - 1},~z_{n})) = \left( \begin{array}{ll}
                                            {\cal O}'& \wp^{n} \\
                                               \wp^{-n + 1} &  {\cal O}'
                                             \end{array} \right) $$

$$\mbox{Stab}(x_{0}) = \left( \begin{array}{ll}
                                               K^{*} & K \\
                                               0     & K^{*}
                       \end{array} \right),
\mbox{Stab}(x_{\infty}) = \left( \begin{array}{ll}
                                               K^{*}  & 0 \\
                                                  K             & K^{*}
                       \end{array} \right).
$$
The stabilizers in the $PGL(V)$ are represented by the matrices from
     $GL(V)$ satisfying the same conditions.
\end{proposition}
\par\smallskip
     
{\sc  Proof  }.  It  is  obvious  that all the vertices from i) belong to
     $\Sigma$. It follows from the theorem 1 (section 1) that  there  are
     no other
vertices. It is also clear that the simplicial complex described in i) is
a minimal complex containing it's vertices.

     The  formulas  for the stabilizers (property ii) can be confirmed by
direct computations.
\par\smallskip
     
     Thus the simplicial structure of an appartment can be  presented  as
the   following  triangulation  of  compactified  line  {\bf  R  }$  \cup
-\infty,~\infty$:
              \\
     \\
\unitlength=1mm
\special{em:linewidth 0.4pt}
\linethickness{0.4pt}
\begin{picture}(132.00,10.00)
\put(10.00,4.00){\circle{4.00}}
\put(20.00,4.00){\circle*{2.00}}
\put(35.00,4.00){\circle*{2.00}}
\put(50.00,4.00){\circle{2.00}}
\put(60.00,4.00){\circle{2.00}}
\put(70.00,4.00){\circle{2.00}}
\put(85.00,4.00){\circle*{2.00}}
\put(100.00,4.00){\circle*{2.00}}
\put(130.00,4.00){\circle{4.00}}
\put(115.00,4.00){\circle{2.00}}
\emline{20.00}{4.00}{1}{35.00}{4.00}{2}
\emline{45.00}{4.00}{3}{49.00}{4.00}{4}
\emline{51.00}{4.00}{5}{59.00}{4.00}{6}
\emline{61.00}{4.00}{7}{69.00}{4.00}{8}
\emline{71.00}{4.00}{9}{75.00}{4.00}{10}
\emline{85.00}{4.00}{11}{100.00}{4.00}{12}
\emline{110.00}{4.00}{13}{114.00}{4.00}{14}
\emline{116.00}{4.00}{15}{120.00}{4.00}{16}
\put(125.00,4.00){\makebox(0,0)[cc]{$\dots$}}
\put(105.00,4.00){\makebox(0,0)[cc]{$\dots$}}
\put(80.00,4.00){\makebox(0,0)[cc]{$\dots$}}
\put(41.00,4.00){\makebox(0,0)[cc]{$\dots$}}
\put(16.00,4.00){\makebox(0,0)[cc]{$\dots$}}
\put(10.00,10.00){\makebox(0,0)[cc]{$x_{0}$}}
\put(20.00,10.00){\makebox(0,0)[cc]{$y_{n - 1}$}}
\put(35.00,10.00){\makebox(0,0)[cc]{$z_{n}$}}
\put(60.00,10.00){\makebox(0,0)[cc]{$x_{i,n}$}}
\put(70.00,10.00){\makebox(0,0)[cc]{$x_{i+1,n}$}}
\put(85.00,10.00){\makebox(0,0)[cc]{$y_{n}$}}
\put(100.00,10.00){\makebox(0,0)[cc]{$z_{n + 1}$}}
\put(115.00,10.00){\makebox(0,0)[cc]{$x_{i,n+1}$}}
\put(130.00,10.00){\makebox(0,0)[cc]{$x_{\infty}$}}
\end{picture}
\\
     \begin{center}
     Pic. 2
     \end{center}
\begin{theorem}
The appartments $\Sigma{.}$ have the following properties:
\begin{quotation}
i) any two simplices are contained in an appartment,

ii) for  any two apppartments
$\Sigma,\Sigma'$  there
exists an isomorphism $i:\Sigma \rightarrow \Sigma'$ such that
$i\mid_{\Sigma \cap \Sigma'} = \mbox{ identity}$,

iii) for any appartment $\bar \Sigma  \subset  \Delta_{.}(G,  K/\bar  K)$
     there exists a unique appartment $\Sigma \subset
\Delta_{.}(G, K/\bar K/k)$ such that $ \pi(\Sigma) = \bar \Sigma, $
 
iv) a geometrical realization $\mid \Sigma_{.} \mid$ of an appartment
$\Sigma_{.}$ is homeomorphic to a closed interval,

v)~$\Sigma_{.} = \{ \sigma \in \Delta_{.} \mid \forall g \in T ~g(\sigma) =
     \sigma \},~N(\Sigma_{.}) \subset \Sigma_{.}$
and the Weyl group  $W$  acts on $\Sigma_{.}$.

If $w \in W$ is an involution then it has a fixed point
$x_{i_{0},n_{0}} \in \Sigma_{0}[2]$ and $w$ is a reflection:
$$w(x_{i,n}) = x_{2i-i_{0},2n-n_{0}},$$
$$w(y_{n_{0}+n}) = z_{n_{0}-n},
~w(z_{n_{0}+n}) = y_{n_{0}-n},~w(x_{0}) = x_{\infty}.$$
If $w \in \Gamma_{K} \cong \mbox{\bf Z} \oplus \mbox{\bf Z} \subset W$ then
$w = (0, 1)$ acts as a shift of the whole structure to the right
$$ w(x_{i,n}) = x_{i,n+2},$$
$$w(y_{n}) = y_{n+2},~w(z_{n}) = z_{n+2},
w(x_{0}) = x_{0},~w(x_{\infty}) = x_{\infty}. $$
The  element  $w  =  (1,  0)$ acts as a shift on the points $x_{i,n}$ but
leaves fixed the points in the inner boundary
$$ w(x_{i,n}) = x_{i+2,n},$$
$$w(y_{n}) = y_{n},~w(z_{n}) = z_{n},
~w(x_{0}) = x_{0},~w(x_{\infty}) = x_{\infty}.  $$
Under the map  $W_{K/{\bar  K}/k}  \rightarrow
W_{K/{\bar K}}$ this action goes to the action of Weyl group  $W_{K/{\bar
K}}$ on an appartment of the tree $\Delta_{.}(K/{\bar K})$.
\end{quotation}
\end{theorem}
\par\smallskip
{\sc  Proof.}  If  we  compare  the  modules  belonging  to an appartment
according  to  proposition  3,  ii)  and  the  modules  belonging  to  an
appartment   of  the  tree  over  a  local  field  of  dimension  1  (see
(\ref{chain}) in section 3) we will see that they will go one to  another
under   the   projection   map.   Thus   we   can   find   an  appartment
$\Sigma$~in~$\Delta_{.}(K/\bar K/k)$,
     projecting onto any given appartment of the tree
$\Delta_{.}(K/\bar K)$.
Note that an appartment always contains an edge from
 $\Delta_{1}[1]$. Then from proposition 2, vi) we get that $\Sigma$
     will be defined in an unique way. We have proved iii).
     
     Let us prove the property i). For any two simplices there  exists  a
subcomplex in $\Delta_{.}(K/\bar K/k)$, having the same combinatorial and
topological  structure  as  the  line  from  picture 2 and containing our
simplices.
 This is obvious from the picture 1. The image of this complex will be an
infinite chain, {\em i. e.} an  appartment  $\bar  \Sigma$  in  the  tree
$\Delta_{.}(K/\bar  K)$  (see section 3). By the result proved above, all
edges of our subcomplex which belong to
$\Delta_{.}[1]$ will also belong to an appartment
$\Sigma$ lying over $\bar \Sigma$. Next we consider the trees which are
the fibers of $\pi$. Looking at them we see that the other  simplices  of
our  subcomplex  also belong to $\Sigma$ (if two appartments in the usual
Bruhat-Tits tree have the same boundary points then they coincide).

     To get ii) we remark that the  intersection  $\Sigma  \cap  \Sigma'$
     will be an "interval", consisting of all the simplices
  lying between two extreme points.
From picture 2 we see the existence of an isomorphism with the properties
which we need.
     
     The  property  iv)  is obvious. In v) we check only the first claim.
The other formulas can be deduced by  direct  computations.  Let  $\sigma
\notin \Sigma_{.}$ and let $\sigma$ be a vertex. Connect this
     vertex  with  $\Sigma$  by  a  minimal  "path"  (=  interval  of  an
appartment). This path will enter into  the  appartment  $\Sigma$  at  an
inner point
      (corollary 1 of proposition 2).
Let  $P$ be a vertex nearest to $\Sigma$ belonging to this path. Then $P$
     belong to a usual Bruhat-Tits tree and there exists
 $g \in G$ such that $g(P) \not = P$. Thus $g(\sigma)
     \not = \sigma$. For the usual tree this property will follow from
the link property (section 3). Namely,  if  $P_{0}$  is  a  point  of  an
appartment  then  the  group  $T$  acts in a simply transitive way on all
edges coming out from
 $P_{0}$ and not lying in the appartment.

     The theorem is proved.
     \par\smallskip
We note that the transformations from the Weyl group  will  be  continous
but not necessarily smooth maps
 of compactified line {\bf R }$ \cup -\infty,
\infty $ into itself.

          If  $P, Q \in \Sigma_{0}$,~then the subcomplex in $\Sigma_{.}$,
containing all the simplices lying between $P$ and $Q$ will be called
 {\em a path} from $P$ to  $Q$ (see picture 2).
\par\smallskip
{\sc Corollary}. {\em For any two vertices $P, Q \in \Delta_{0},~P \not =
Q$ there exists a unique path $PQ$ between them}.
\par\smallskip
As in the usual theory of the Bruhat-Tits  tree  we  can  introduce  some
intrinsically defined metric over our tree.

If  $<L>,<L'> \in \Delta_{0}(G)[2]$ then by Cartan decomposition (theorem
2 of section 2) there exists a basis
$e_{1},  e_{2}$ in $V$ such that
$$L = {\cal O}'e_{1} \oplus {\cal O}'e_{2},
~L' = a_{1}e_{1} \oplus  a_{2}e_{2}, $$
where  $a_{1},  a_{2}$  are  some  fractional  ${\cal   O}'$-ideals   and
$\nu'(a_{1}) \geq \nu'(a_{2})$.
\par\smallskip
 {\sc Definition 8}
~$d( <L>, <L'> ) = \nu'(a_{1}) - \nu'(a_{2}), $
where $<L>,<L'>$ are two vertices from $\Delta_{0}[2]$.
\par\smallskip
\begin{theorem}.
The function $ d(.,.)$ is a correctly defined metric on the set
$\Delta_{0}[2]$ having non-archimedean values in $\Gamma_{K}^{+}$. It
has the following properties
\begin{quotation}
i) $d(.,.)$ is invariant under the action of $G$.

ii)the projection map $\pi$ is a distance-decreasing map, precisely
$$ d(\pi(x), \pi(y)) = \pi (d(x, y))$$

iii)for any appartment $\Sigma$ there exists a simplicial map
$ \rho: \Delta_{.} \rightarrow \Sigma_{.} $
which  is  a  retraction onto $\Sigma$ and which is a distance-decreasing
map on subset $\Delta_{.}[2]$.
     
iv)  let  $u,  t$  be  local  parameters  of  the field $K$ and $P, Q \in
\Delta_{0}[2]$. Then $d(P, Q) = (m, n)$ and we have
     $$ n = d(\pi P, \pi Q) \mbox{~in~} \Delta(K/\bar K)[1],$$
     $$ m = d(Q, Q') \mbox{~in~}
\Delta(\bar K/k)[1] \cong
\pi^{-1}(Q) \cap \Delta(K/\bar K/k)[2],  $$
     where $Q' =
 \left( \begin{array}{ll}
                          t^{n}  & 0\\
                           0 & 1
                           \end{array} \right) P \in \pi^{-1}(Q) $

v) if $R \in PQ$, then $d(P, R) + d(R, Q) = d(P, Q)$,

vi) for $P, Q, P', Q' \in \Delta_{0}[2]$ there exist $g \in G$ such that
     $gP = P', gQ = Q'$~if and only if~$d(P, Q) = d(P', Q')$.
\end{quotation}
\end{theorem}
\par\smallskip
{\sc Proof.} If we change a module inside it's class then the same number
will be added to the $\nu'(a_{1})$ and $\nu'(a_{2})$. Consequently, their
difference will be unchanged. The properties i) and ii) follows directly
     from the definition. Let us show how to construct the retracting map.

     Take an edge $\sigma$ of the appartment $\Sigma$. Then
$\Sigma_{.} - \sigma$ can be decomposed into two pieces
      $\Sigma_{+}$ and $\Sigma_{-}$.
Let $0$ and $\infty$ be the points from  the  external  boundary  of  the
appartment. We assume that $0~(\infty)$ are the limit points for
 $\Sigma_{+}~(\Sigma_{-})$.

For  any  point  $P  \in  \Delta_{0}[0]$  which  does  not  belong to the
appartment there is a unique shortest path which connects $P$  with  some
point
$Q(P)$ of the appartment.

Thus the whole external boundary $\Delta_{.}[0]$  can be divided into two pieces
$\Delta[0]_{+}$ and $\Delta[0]_{-}$.
The  first  piece $\Delta[0]_{+}$ will contain 0 and all the points which
 are connected with $\Sigma_{+}$. All the other points will belong to
 $\Delta[0]_{-}$. We start to construct $\rho$ from the external boundary:
 $$ \rho (\Delta[0]_{+}) = 0, ~\rho (\Delta[0]_{-}) = \infty .$$
     Then if $P \in \Delta[0]_{+}, ~P \not  =  0$  there  are  two  paths
connecting the point $Q(P)$ with external boundary: the path between
$P$~and~$Q(P)$, and the path between $0$~and~$Q(P)$ (a part of $\Sigma_{+}$).
There  exists  a  unique  simplicial  bijection  $s_{P}$ of one path onto
another one. Let us put

     $$\rho (\sigma) = s_{P}(\sigma), \mbox{~if}~\sigma \mbox{~lies on the
path between }P~\mbox{and}~Q(P).$$
     The same definition works for $\Delta[0]_{-}$. It is straightforward
that the constructed map is correctly  defined  on  the  whole  tree  and
satisfies all the conditions from iii).
     
     Properties  iv)  and  v)  follows  from  ii) and direct computations
(compare with
     theorem 5, v) ). To get vi) we first observe that we can assume
     $P' = P$ (since $G$ is transitive on the tree) and $\pi(Q) = \pi(Q')$
  (apply the same property for the tree $\Delta(K/\bar K)$). Now let
     $n = d(P, Q) = d(P, Q')$~and we assume $n > 1$(otherwise we are done
by the property v) for the tree $\Delta(\bar K/k)$).Let $R$ be  a  common
point of the paths
     $PQ$~and~$PQ'$ such that the intersection of the paths $RQ$ ~and~$RQ'$
     is  $R$.  Let  us  denote by $R_{0}$ the inner boundary point of the
path $RP$ which is closest to $R$. Then $R, R_{0}$  belong  to  the  same
fiber as $Q$~and~$Q'$ and the
     equality $d(P, Q) = d(P, Q')$ is equivalent to $d(R, Q) = d(R, Q') $
     in the tree $\Delta(\bar K/k)[1] \cong \Delta' =
\pi^{-1}(\pi (Q))  \cap  \Delta(K/\bar  K/k)[2]$.  By  proposition  3,  ii) the
stabilizer $G'$ of the points $R_{0}$~and $P$ has the matrix form
$ \left( \begin{array}{ll}
                          {\cal O'}  & {\cal O}\\
                           m_{i, n} &    {\cal O'}
                   \end{array} \right)  $~for some $i$. By the corollary 2
     of  proposition  2  $G'$  acts  on  $\Delta'$  as  a  group of upper
     triangular matrices
$ \left( \begin{array}{ll}
                          {\cal O}^{*}_{\bar K}  & {\bar K}\\
                           0  &    {\cal O}^{*}_{\bar K}
                   \end{array} \right)  $. This group acts
     transitively on the boundary of $\Delta'$ outside $R_{0}$  and  thus
under our distance condition it will move $Q$ to $Q'$.
     
     The theorem is proved.
                           
\par\smallskip
The  last  general notion which will be mentioned here is the type of the
vertices and also of the simplices. Let us consider an exact sequence
$$0 \rightarrow PGL^{+}(V) \rightarrow PGL(V) \rightarrow
\Gamma_{K}/2\Gamma_{K} \rightarrow 0$$
where the right hand map is $\nu'(det(.))$ {\em mod} 2. As we know
$$\Gamma_{K}/2\Gamma_{K} \cong \mbox{\bf Z}/2\mbox{\bf Z} \oplus
\mbox{\bf Z}/2\mbox{\bf Z}$$
It can be shown that the stabilizers of the vertices belong to the subgroup
$PGL^{+}(V)$ and thus we have a canonical map
$$ \Delta_{0}[2] \rightarrow \Gamma_{K}/2 \Gamma_{K}$$
which assign to the vertices {\em four} possible values, their type.  The
type of a simplex will be then a subset of $\Gamma_{K}/2 \Gamma_{K}$. The
type  is invariant under the action of $SL(V)$ and the fundamental domain
of this action is a disjoint union of two edges which are mapped  by  the
projection  map  on  the  adjacent  vertices  of  an  edge  in  the  tree
$\Delta_{.}(K/{\bar K})$.
     
The integer points of the lattice $\Gamma_{K}$ can be located on an  real
plane and it seems more reasonable to have a srtucture of dimension 2 on
     our simplicial set.

In the case of local field $K$ of dimension $n$
 the number of types equals to $2^{n}$
and the building of $G$ (see\cite{P}) could have a dimension depending on
 $n$.

But  if  we  would like to preserve one of most important features of the
Tits theory - the geometrical structure of reflections,  walls,  chambers
and  so on then we are forced to introduce the simplicial structure as we
did above. The reason is that the involutions from the  Weyl  group  have
very  small  fixed  point  set on the lattice $\Gamma_{K}$ (see theorem 3
above). In particularly, in dimension two they have  the  points  as  the
fixed  points  but not the lines as would be the case if the dimension of
our building were two.

{\bf Remark 6}.
Our use of the topology was rather  artificial.  It  seems  there  should
exist    a    purely    simplicial    construction    which   binds   the
$\Delta_{.}[m]$-pieces  of  the  tree   together.   We   can   define   $
\Delta_{.}(G) = \Delta_{.}[2] * \Delta_{.}[1]
* \Delta_{.}[0], $
where * is a join of the simplicial complexes. Then the group
$G \times G \times G $
will act on the whole $\Delta_{\mbox{max}}(G)$ in a transitive way and we
will  have  a  one to one correspondence between subgroups of this larger
group and the simplexes of the new complex which has a dimension  4.  The
same  remark  is  true for the groups of higher rank over arbitrary local
fields \cite{P}.
                                                               
We also add the following problem.

{\sc Problem 3.} It is  well  known  that  the  buildings  of  the  group
$\mbox{PGL}(V)$ (and in particularly the Bruhat-Tits tree) can be defined
as  classes  of norms on the space $V$ \cite[II, 1.1]{BT2, S2} . There is
no doubt that  this  approach  can  be  developed  also  for  the  higher
buildings of this group also. But this should give directly a geometrical
realization  of  the  simplicial set $\Delta_{.}(G)$ which was defined in
\cite{P}.


\begin{thebibliography}{99}


\bibitem[\bf 1]{B}
 Bourbaki N., {\em Groupes  et  Algebre  de  Lie}  (chapitre  IV  -  VI),
   Hermann, Paris, 1968


\bibitem[\bf 2]{BT1}
 Bruhat  F.,  Tits  J.,  {\em Groupes r\'eductives sur un corps local. I.
      Donn\'ees radicielles valu\'ees}, Publ. Math. IHES, {\bf 41}(1972),
      5-251;{\em
   II.Sch\'emas  en  groupes,   Existence   d'une   donn\'ee   radicielle
        valu\'ee}, Publ. Math. IHES, {\bf 60}(1984), 5-184

\bibitem[\bf 3]{BT2}
  Bruhat  F., Tits J., {\em Sch\'emas en groupes et immeubles des groupes
                  classiques  sur  un  corps  local},  Bull.  Soc.  Math.
                  France. {\bf 112}(1984), 259-301

\bibitem[\bf 4]{C1}   Cartier P., {\em G\'eom\'etrie et Analyse sur les Arbres},
S\'eminaire Bourbaki, 1971/1972,
            n   407,   Lecture   Notes   in   Mathematics,   {\bf   317},
            Springer-Verlag, Berlin, 123-140

\bibitem[\bf 5]{C2}   Cartier P., {\em Fonctions harmoniques sur un arbre},
Symposia
                    Matematica, {\bf 9} (1972), 203-270

\bibitem[\bf 6]{D}   Deligne P., {\em The\'orie de Hodge II}., Publ. Math. IHES,
{\bf 40}(1971), 5-58

\bibitem[\bf 7]{FP}
Fimmel T., Parshin A. N., {\em An Introduction  into  the  Higher  Adelic
Theory} (in preparation)
     
\bibitem[\bf 8]{G}
 Gantmacher F.R., {\em Theory of matrices}, Nauka, Moscow, 1966

\bibitem[\bf 9]{GI}
 Goldman O., Iwahori N., {\em The space of $\wp$-adic norms}, Acta Math.,
      {\bf 109}(1963), 137-177

\bibitem[\bf 10]{IM}
     Iwahori  N., Matsumoto H., {\em On some Bruhat decomposition and the
     structure of the Hecke rings of the  $\wp$-adic  Chevalley  groups},
     Publ. Math. IHES, {\bf 25}(1965), 5-48
     

\bibitem[\bf 11]{P}
 Parshin A. N. {\em Higher Bruhat-Tits buildings and Vector Bundles on an
                       Algebraic Surface}, Algebra and Number Theory
        (Proc.Conf.  held  at  the Inst.Exp.Math. , Univ.Essen, Dec. 2-4,
                       1992), de Gruyter, Berlin, 1994, 165-192



\bibitem[\bf 12]{R}
 Ronan M., {\em Buildings: Main Ideas and Applications.I.}, Bull.  London
          Math.Soc. {\bf 24}(1992), 1-51;{\em II.},
        Bull. London Math.Soc.
          {\bf 24}(1992), 97-126


\bibitem[\bf 13]{S1}
      Serre J.-P., {\em Corps locaux}, Hermann, 1968
     
\bibitem[\bf 14]{S2}
  Serre  J.-P.,  {\em Arbres, amalgames}, $SL_{2}$, Asterisque, {\bf 46},
          Soc. Math. France, Paris, 1977

\bibitem[\bf 15]{T}
 Tits J., {\em On Buildings and their applications}, Proc. Intern. Congr.
             Math. (Vancouver 1974), Canad. Math. Congr., Montreal, 1975,
        Vol. {\bf 1}, 209-220

\end{thebibliography}
\end{document}